\documentstyle[12pt,amstex,amssymb]{article}
\newtheorem{lemma}{Lemma}

\newcommand{\average}[1]{\left\langle #1 \right\rangle}
\newcommand{\brk}[1]{\left( #1 \right)}
\newcommand{\Brk}[1]{\left[ #1 \right]}
\newcommand{\BRK}[1]{\left\{ #1 \right\}}
\newcommand{\scalar}[2]{\left( #1 , #2 \right)}
\newcommand{\cov}[2]{\text{Cov}\left[ #1 , #2 \right]}
\newcommand{\half}{\frac{1}{2}}
\newcommand{{\V}}{\mathcal{M}}
\numberwithin{equation}{section}

\begin{document}
\begin{center}

{\Large\bf Prediction of Large-Scale Dynamics}\medskip\\
{\Large\bf Using Unresolved Computations}

\bigskip
{\sc Alexandre J. Chorin, Anton Kast, and Raz Kupferman}\footnote{
Current addrress: Institute of Mathematics,
The Hebrew University,
Jerusalem  91904 Israel.}\bigskip\\
Department of Mathematics \\
Lawrence Berkeley National Laboratory \\
Mail Stop 50A-2152 \\
1 Cyclotron Road \\
Berkeley, CA 94720\\
{\tt chorin@@math.berkeley.edu}\\
{\tt anton@@math.lbl.gov}\\
{\tt raz@@math.lbl.gov}
\end{center}

\bigskip\bigskip

{\bf Abstract.}
We present a theoretical framework and numerical methods for
predicting the large-scale properties of solutions of partial
differential equations that are too complex to be properly
resolved. We assume that prior statistical information about the
distribution of the solutions is available, as is often the case
in
practice. The quantities we can compute condition the prior
information and allow us to calculate mean properties of solutions
in
the future.  We derive approximate ways for computing the
evolution
of the probabilities conditioned by what we can compute, and
obtain
ordinary differential equations for the expected values of a set
of
large-scale variables.  Our methods are demonstrated on two
simple but
instructive examples, where the prior information consists of
invariant canonical distributions

\vfill
\noindent This work was supported in part by the US Department of
Energy under
contract\linebreak DE--AC03--76--SF00098, and in part by the National Science
Foundation under grant\linebreak DMS 94--14631.

\begin{center}Subject Classification: Primary 65M99\end{center}

\newpage

\section{Introduction}

There are many problems in science that can be modeled by a set of
differential equations, but where the solution of these equations
is
so complicated that it cannot be found in practice, either
analytically or numerically.  For a numerical computation to
be
accurate the problem must be well resolved, i.e, enough variables
(or
``degrees of freedom'') must be represented in the calculation to
capture all the relevant features of the solution; insufficient
resolution yields sometimes disastrous results.  A well-known
example
in which good resolution cannot be achieved is turbulent flow,
where
one has to resolve all scales ranging from the size of the system
down
to the dissipation scale---a prohibitively expensive requirement.
One is then compelled to consider the question of how to predict
complex behavior when the number of variables that can be used in
the
computation is significantly less than needed for full resolution.
This is the question considered in the present paper; part of the
theoretical framework and methods have already been briefly
discussed
in \cite{CKK98}.

Studies on underresolved problems exist in a wide range of
different
contexts, along with a large amount of literature that describes
problem-specific methods.  In turbulence, for example, there
are
various modeling methods for large eddy simulations.  In all
cases one
needs to make additional assumptions about the relation between
those
degrees of freedom that are represented in the computation and the
``hidden'', or ``invisible'' degrees of freedom that are discarded
from the computation.  A number of interesting attempts have
been made
over the years to fill in data from coarse grids in difficult
computations so as to enhance accuracy without refining the grid
(see
e.g.~\cite{SM97,MW90}).  Indeed, nothing can be done without
some
information regarding the unresolved degrees of freedom.
Such
additional assumptions are usually motivated by intuitive
reasoning and their validity is usually assessed by comparing the
resulting predictions to experimental measurements.

In many problems the lack of resolution is due primarily to the
insufficiency and sometimes also the inaccuracy of the
measurements
that provide initial conditions for the system of equations.
This is
the case for example in weather forecasting, where the initial
information consists of local weather measurements collected at a
relatively small number of meteorological stations.  The problem
of
insufficient and sometimes noisy data is not considered in the
present
paper.  We focus here on the case where underresolution is
imposed by
computational limitations.  Initial data will be assumed to
be
available at will, and this assumption will be fully exploited by
allowing us to select the set of degrees of freedom that are
represented in the computation at our convenience.  Another
issue that
often arises in the modeling of complex systems is uncertainty
regarding the equations themselves.  This important question is
also
beyond the scope of this paper; the adequacy of the system of
equations to be solved is taken for granted.

We now define the problem and introduce some of the nomenclature:
We
consider a system described by a differential equation of the form
\begin{equation}
u_t = F(u),
\label{TheGeneralForm}
\end{equation}
where $t$ is time, subscripts denote differentiation, $u(x,t)$ is
the
dependent variable, and $F(u)= F(u,u_x,u_{xx},\ldots)$ is a
(generally
nonlinear) function of its arguments; the spatial coordinate $x$
and
the dependent variable $u$ can be of arbitrary dimensionality.

To solve an equation of the form (\ref{TheGeneralForm}) on a
computer
one ordinarily discretizes the dependent variable $u(x,t)$ both in
space and time and replaces the differential equation by an
appropriate relation between the discrete variables.  As
described,
the solution to the discrete system may approximate the solution
of
the differential equation well only if the discretization is
sufficiently refined.  It is our basic assumption that we
cannot
afford such a refined discretization, and must therefore be
content
with a much smaller number of variables.  One still has the
liberty to
choose the degrees of freedom that are retained in the
computation;
those will be chosen, for convenience, to be linear functionals of
the
dependent variable $u(x,t)$:
\begin{equation}
U_\alpha[u(\cdot,t)] \equiv
\scalar{g_\alpha(\cdot)}{u(\cdot,t)} \equiv
\int g_\alpha(x) u(x,t) \, dx,
\label{CollectiveVariables}
\end{equation}
where $\alpha$ is an index that enumerates the selected degrees of
freedom.  Variables of the form (\ref{CollectiveVariables}) will
be
referred to as \emph{collective variables}; every collective
variable
$U_\alpha$ is defined by a \emph{kernel} $g_\alpha(x)$.  Point
values
of $u(x)$ at a set of points $x_\alpha$, and spectral components
of
$u(x)$ for a set of modes $k_\alpha$ are two special cases of
collective variables; in the first case the corresponding kernels
are
delta functions, $g_\alpha(x) = \delta(x-x_\alpha)$, whereas in
the
second case the kernels are spectral basis functions,
$\exp(ik_\alpha\cdot x)$.  We assume that our computational
budget
allows us to operate on a set of at most $N$ collective variables,
so
that $\alpha = 1,\ldots,N$.  The question is, what can be
predicted
about the state of the system at a future time $t$ given the values
of
the collective variables $U_\alpha$ at an initial time $t=0$?

Suppose that we know at time $t=0$ that the collective variables
$U_\alpha$ assume a set of values $V_\alpha$. (We will denote by $U
=
(U_1,\ldots,U_N)^T$ and $V = (V_1,\ldots,V_N)^T$ the vectors whose
entries are the collective variables and their initial values,
respectively.)  Our postulate that the number of collective
variables
$N$ does not suffice to resolve the state of the system implies
that
the initial data, $V$, do not determine sharply enough the initial
condition, $u(x,0)$. A priori, every function $u(x,0)$ that is
compatible with the given values of the collective variables, that
is,
belongs to the set
\begin{equation}
{{\V}(V)} =
\left\{ v(x) : U_\alpha[v(\cdot)]  = V_\alpha,
\quad \alpha=1,\ldots,N \right\}.
\label{SetOfIC}
\end{equation}
is a plausible initial condition. One could define underresolution
in
terms of the set of functions (\ref{SetOfIC}); the problem is
underresolved if this set is non-trivial.  Clearly, the state of
the
system at future times depends on the particular initial condition;
in
many cases it is even very sensitive to small variations in the
initial condition. One wonders then in what sense the future can
be
predicted when the initial condition is not known with certainty.

The essence of our approach is the recognition that
underresolution
necessarily forces one to consider the evolution of a set, or
ensemble, of solutions, rather than a single initial value
problem.
This requires the replacement of equation (\ref{TheGeneralForm}) by
a
corresponding equation for a probability measure defined on the
space
of the solutions of (\ref{TheGeneralForm}).  The prediction of
the
future state of the system can then be reinterpreted as the
prediction
of most likely, or mean, properties of the system.  Loosely
stated, in
cases where sufficient resolution cannot be achieved the original
task
of solving an initial value problem has to be replaced by a more
modest one---the determination of ``what is most likely to happen
given what is initially known.''

At first, there seems to be no practical progress in the above
restatement of the problem.  First, the statistical problem
also
requires initial conditions; a measure defined on the space of
initial
conditions $u(x,0)$ must be provided for the statistical problem to
be
well-defined.  Second, the high-dimensional Liouville equation
that
describes the flow induced by (\ref{TheGeneralForm}) is not easier
to
solve than the original initial value problem.  It turns out
that in
many problems of interest there exists a natural measure $\mu$
that
characterizes the statistical properties of the system; what is
meant
by ``natural'' has to be clarified; an important class of such
measures are invariant ones.  We are going to use this
information to
partially cure the two aforementioned difficulties: First, this
measure will define the initial statistical state of the system by
being interpreted as a ``prior'' measure---a quantification of our
beliefs regarding the state of the system prior to the
specification
of any initial condition.  The initial values of the
collective
variables are constraints on the set of initial states and induce
on
$\mu$ a conditional measure that constitutes an initial condition
for
the Liouville equation.  Second, the existence of a
distinguished
statistical measure suggests a way to generate a hierarchy of
approximations to the Liouville equation, examples of which will
be
described in the following sections.

The rest of this paper is organized as follows: In Section~2 we
present our theory, and provide a recipe
(\ref{TheEffectiveEquation})
for approximating the mean evolution of a set of collective
variables. In Section~3 we derive formulas for the calculation of
conditional expectations in the case of Gaussian prior measures;
these
are necessary for the evaluation of the right-hand side of
equation
(\ref{TheEffectiveEquation}). In Sections 4 and 5 we demonstrate
the
power of our theory by considering two examples: a linear
Schr\"odinger equation and a nonlinear Hamiltonian system.
Conclusions
are presented in Section~6.

\section{Presentation of the theory}

Our starting point is a general equation of motion of the form
(\ref{TheGeneralForm}), and a set of collective variable
$U_\alpha$
defined by (\ref{CollectiveVariables}) for a set of kernels
$g_\alpha(x)$; the question of what constitutes a good choice of
kernels will be discussed  below.

In many problems of interest there exists a measure on the space
of
solutions of (\ref{TheGeneralForm}) that is invariant under the
flow
induced by (\ref{TheGeneralForm}); a measure that has this property
is
referred to as an \emph{invariant measure}. Invariant measures are
known to play a central role in many problems; macroscopic systems
(that is, systems that have a very large number of degrees of
freedom)
whose macroscopic properties do no change in time, often exhibit
an
invariant statistical state.  By that we mean the following:
when the
large scale observable properties of the system remain constant in
time, the likelihood of the microscopic degrees of freedom to be
in
any particular state is distributed according to a measure that is
invariant in time.  We will assume that such an invariant
measure
$\mu_0$ exists and that we know what it is. The measure $\mu_0$
will
then be postulated to be the \emph{prior measure}, i.e, it
describes
the probability distribution of initial conditions before any
measurement has been performed.  We will denote averages with
respect
to the invariant measure $\mu_0$ by angle brackets
$\average{\cdot}$;
let $O[u(\cdot)]$ be a general functional of $u$, then
\begin{equation}
\average{O} = \int O[u(\cdot)] \,d\mu_0,
\end{equation}
where the integration is over an appropriate function space.
We shall write formally,
\begin{equation}
d\mu_0 = f_0[u(\cdot)] \,[du],
\end{equation}
as if the measure $\mu$ were absolutely continuous with respect to
a
Lebesgue measure, where $f_0[u]$ is the invariant probability
density,
and $[du]$ is a formal product of differentials.

We next assume that a set of measurements has been carried out and
has
revealed the values $V_\alpha$ of the collective variables
$U_\alpha$
at time $t=0$. This information can be viewed as a set of
constraints
on the set of initial conditions, which is now given by
(\ref{SetOfIC}). Constraints on the set of functions $u(x)$
automatically induce on $\mu_0$ a \emph{conditional measure}, which
we
denote by $\mu_V$. In a physicist's notation,
\begin{equation}
d\mu_V = f_V[u(\cdot)]\,[du] =
c f_0[u(\cdot)] \, [du] \times \prod_{\alpha=1}^N
\delta \brk{ U_\alpha[u(\cdot)] - V_\alpha},
\label{InitialDensity}
\end{equation}
where $f_V[u(\cdot)]$ is the conditional probability density, and
$c$
is an appropriate normalization factor. The conditional
probability
density is equal, up to a normalization, to the prior probability
density projected on the space of functions ${\V}(V)$ that are
compatible with the initial data. Note that the conditional
measure
$\mu_V$ is, in general, not invariant.  Averages with respect to
the
conditional measure will be denoted by angle brackets with a
subscript
that symbolizes the constraints imposed on the set of functions,
\begin{equation}
\average{O}_V \equiv \int O[u(\cdot)] f_V[u(\cdot)] \,[du].
\end{equation}

The dynamics have not been taken into consideration so far, except
for
the fact that the measure $\mu_0$ was postulated to be
invariant.  Let
$f[u(\cdot),t]$ be the probability density of the solutions of
(\ref{TheGeneralForm}) at time $t$, that is, the probability
density
that evolves from the initial probability density $f_V[u(\cdot)]$
under the flow induced by (\ref{TheGeneralForm}); it satisfies the
Liouville equation \cite{Ris84}
\begin{equation}
f_t + \scalar{\frac{\delta f}{\delta u}(\cdot)}{F(u(\cdot))} = 0,
\label{LiouvilleEquation}
\end{equation}
where $\frac{\delta f}{\delta u}$ denotes a functional derivative.
An
equivalent statement is that if $S_t$ denotes the time evolution
operator induced by (\ref{TheGeneralForm}), i.e., $S_t : u(x,0)
\to
u(x,t)$, then
\begin{equation}
f[u(\cdot),t] = f[S^{-1}_t u(\cdot), 0] = f_V[S^{-1}_t u(\cdot)],
\label{DensityAtTimeT}
\end{equation}
where $S^{-1}_t$ is the operator inverse to $S_t$, which we assume
to
exist.

The objective that has been defined in the introductory section is
to
calculate the expectation value of observables $O[u(\cdot)]$ at
time
$t$, given the initial data $V$. In terms of the notations
introduced
above this is given by
\begin{equation}
\average{O[u(\cdot),t]}_V = \average{O[S_t u(\cdot)]}_V
\label{Def1}
\end{equation}
(operators are generally treated as function of the dependent
variable
and time, $O[u(\cdot),t]$; when no reference to time is being made
the
expression refers to the initial time).

We next make the following observations: (i) The initial
probability
measure (\ref{InitialDensity}) is completely determined by the $N$
numbers $V_\alpha$. (ii) By the invariance of $f_0[u]$ and by
equation
(\ref{DensityAtTimeT}), the probability density at later time $t$
can
still be represented as the invariant density projected on a set
of
$N$ conditions; specifically,
\begin{equation}
f[u(\cdot),t] = c \, f_0[u(\cdot)] \prod_{\alpha=1}^N
\delta\Brk{ \scalar{g_\alpha(\cdot)}{S^{-1}_t u(\cdot)} -
V_\alpha}.
\label{DensityAtTimeT2}
\end{equation}
Note however that the set of functions that support this measure
at
time $t$ is generally not of the form (\ref{SetOfIC}), that is,
the
observable $\scalar{g_\alpha(\cdot)}{S^{-1}_t u(\cdot)}$ is not a
linear functional of $u$.

These observations suggest an approximate procedure for solving
the
Liouville equation (\ref{LiouvilleEquation}). We propose an ansatz
in
which the $N$ conditions that are imposed on $\mu_0$ remain for
all
times conditions on the values of the collective variables $U$;
namely, the probability density is specified by a time-dependent
vector of $N$ numbers $V_\alpha(t)$, such that
\begin{equation}
f[u(\cdot),t] \approx  c\,\, f_0[u(\cdot)] \prod_{\alpha=1}^N
\delta\Brk{ U_\alpha[u(\cdot)] - V_\alpha(t)}.
\label{TheAnsatz}
\end{equation}

One has still to specify the time evolution of the vector $V(t)$.
Suppose that the distribution of solutions is indeed given by
(\ref{TheAnsatz}) at time $t$, and consider a later time $t+\Delta
t$. The value of the observable $U_\alpha[u(\cdot)]$ at the later
time
will, in general, not be uniform throughout the ensemble of
solutions. The ansatz (\ref{TheAnsatz}) projects the distribution
back
onto a set of solutions ${\V}(V(t+\Delta t))$. A natural choice
for $V_\alpha(t+\Delta t)$ is the expectation value of the
collective
variable $U_\alpha[u(\cdot)]$ given that the distribution at time
$t$
was (\ref{TheAnsatz}):
\begin{equation}
\begin{split}
V_\alpha(t+\Delta t) & \approx
\average{U_\alpha[u(\cdot),t+\Delta t]}_{V(t)} =
\\ &=
\average{U_\alpha[u(\cdot)]}_{V(t)} + \Delta t\,
\average{\scalar{g_\alpha(\cdot)}{F(u(ot))}}_{V(t)} + O(\Delta t^2).
\end{split}
\end{equation}
Taking the limit $\Delta t\to0$ we finally obtain,
\begin{equation}
\frac{d V_\alpha}{d t} =
\average{\scalar{g_\alpha(\cdot)}{F(u(ot))}}_{V(t)}.
\label{TheEffectiveEquation}
\end{equation}

Equation (\ref{TheEffectiveEquation}) is our main tool in the
present
paper and we will next discuss its implications:

\begin{itemize}

\item
Equation (\ref{TheEffectiveEquation}) constitutes a closed set of
$N$
ordinary differential equations, which by our postulate is within
the
acceptable computational budget.

\item
The central hypothesis in the course of the derivation was that
the
distribution of solutions can be approximated by
(\ref{TheAnsatz}).
This approximation assumes that for all times $t$ the collective
variable $U_\alpha$ has a uniform value $V_\alpha$
\emph{for all the trajectories in the ensemble of solutions}.
This
assertion is initially correct (by construction) at time $t=0$,
but
will generally not remain true for later times.  The
approximation is
likely to be a good one as long as the above assertion is
approximately true, that is, as long as the distribution of values
assumed by the collective variables remains sufficiently
narrow.  In
many cases it is possible to guarantee a small variance by a
clever
selection of collective variables (i.e., of kernels). Note
furthermore
that the smallness of the variance can be verified
self-consistently
from the knowledge of the probability density (\ref{TheAnsatz}).

\item
Equation (\ref{TheEffectiveEquation}) still poses the technical
problem of computing its right-hand side.  This issue is the
subject
of the next section.

\item
The case where the equations of motion (\ref{TheGeneralForm}) are
linear, i.e,
\begin{equation}
u_t = L u,
\end{equation}
with $L$ being a linear operator, can be worked out in detail.
Using
the fact that $S_t = \exp(L t)$, the solution to the Liouville
equation (\ref{DensityAtTimeT2}) can be rearranged as
\begin{equation}
f[u(\cdot),t] = c f_0[u(\cdot)] \,\prod_{\alpha=1}^N
\delta\Brk{\scalar{e^{-L^\dagger t} g_\alpha(\cdot)}{u(\cdot)} -
V_\alpha},
\end{equation}
where $L^\dagger$ is the linear operator adjoint to $L$.  Thus,
the
probability density for all times is $f_0$ projected on the set of
functions for which a set of $N$ linear functionals of $u$ have
the
values $V$; note that $V$ here is not time dependent, but is the
vector of initial values of the collective variables $U$.  The
kernels
that define these functionals are time dependent, and evolve
according
to the dual equation
\begin{equation}
\frac{d g_\alpha}{d t} = - L^\dagger  g_\alpha.
\end{equation}
If the kernels $g_\alpha$ are furthermore eigenfunctions of the
dual
operator $L^\dagger$ with eigenvalues $\lambda_\alpha$, the ansatz
(\ref{TheAnsatz}) is exact, with $V_\alpha(t) =
V_\alpha(0)\,e^{\lambda_\alpha t}$.  Hald \cite{Hal98} shows
that by
selecting kernels that are \emph{approximate} eigenfunctions of
$L^\dagger$, one can bound the error introduced by the ansatz
(\ref{TheAnsatz}), while retaining the simplicity of the
procedure.

\item
The two alternatives of evolving either the values $V_\alpha$ or
the
kernels $g_\alpha(x)$ are analogous to Eulerian versus Lagrangian
approaches in fluid mechanics, or Schr\"odinger versus
Heisenberg
approaches in quantum mechanics.  For nonlinear equations one
has a
whole range of intermediate possibilities; for example one may
split
the operator $F$ in equation (\ref{TheGeneralForm}) as $F = L +
Q$,
where $L$ is linear.  The kernels can be evolved according to
the
linear operator, while the values of the collective variables can
be
updated by the remaining nonlinear operator.  The art is to
find
partitions $F = L + Q$ that minimize the variance of the
distribution
of values assumed by the collective variables.

\item
Equation (\ref{TheEffectiveEquation}) should be viewed as a first
approximation to the solution of the Liouville equation, where the
only information that is updated in time is the mean value of a
fixed
set of collective variables. In principle, one could also update
higher moments of those variables, and use this additional
information
to construct a better approximation. For example, equipped with
the
knowledge of means and covariances one could find new kernels and
new
values for the corresponding collective variables, such that the
distribution obtained by conditioning the invariant distribution
with
those new constraints is compatible with the calculated means and
covariances. Thus, one could imagine an entire hierarchy of
schemes
that take into account an increasing number of moments of the
resolved
variables.
\end{itemize}

\section{Conditional expectation with Gaussian prior}

Equation (\ref{TheEffectiveEquation}) is a closed set of equations
for
the vector $V(t)$, which requires the computation of a conditional
average on its right-hand side. To have a fully constructive
procedure, we need to evaluate conditional averages
$\average{O[u(\cdot)]}_V$, where $O$ is an arbitrary observable,
and
$V$ denotes as before the vector of values of a set of collective
variables $U$ of the form (\ref{CollectiveVariables}). In this
section
we present three lemmas that solve this problem for the case where
the
prior measure $\mu_0$ is Gaussian. In the two examples below, the
prior measure is either Gaussian or can be viewed as a perturbation
of
a Gaussian measure.

The random function $u(x)$ has a Gaussian distribution if its
probability density is of the form
\begin{equation}
f_0[u(\cdot)] = Z^{-1} \exp\brk{
-\half \iint u(x) a(x,y) u(y)\,dx\,dy
+ \int b(x) u(x) \, dx},
\end{equation}
where $a(x,y)$ and $b(x)$ are (generalized) functions, and $Z$ is
a
normalizing constant. The functions $a(x,y)$ and $b(x)$ are related
to
the mean and the covariance of $u(x)$ by
\begin{equation}
\average{u(x)} = \scalar{a^{-1}(x,\cdot)}{b(\cdot)},
\end{equation}
and
\begin{equation}
\cov{u(x)}{u(y)} \equiv
\average{u(x) u(y)} - \average{u(x)}\average{u(y)} = a^{-1}(x,y),
\label{UncCovariance}
\end{equation}
where the generalized function $a^{-1}(x,y)$ is defined by the
integral relation
\begin{equation}
\scalar{a(x,\cdot)}{a^{-1}(\cdot,y)} =
\scalar{a^{-1}(x,\cdot)}{a(\cdot,y)} = \delta(x-y).
\end{equation}

To compute the expectation value of higher moments of $u$ one can
use
Wick's theorem \cite{Kle89}:
\begin{multline}
\average{(u_{i_1}-\average{u_{i_1}}) \cdots
(u_{i_l}-\average{u_{i_l}})} = \\
\left\{
\begin{aligned}
0, & \qquad l \text{ odd} \\
\sum
\cov{u_{i_{p_1}}}{u_{i_{p_2}}} \cdots
\cov{u_{i_{p_{l-1}}}}{u_{i_{p_l}}}, & \qquad l \text{ even}
\end{aligned}
\right. ,
\label{Wick}
\end{multline}
with summation over all possible pairings of $\{i_1,\ldots,i_l\}$.

Next, suppose that the random function $u(x)$ is drawn from a
Gaussian distribution, and a set of measurements reveal the vector
of
values $V$ for a set of collective variables $U$ of the form
(\ref{CollectiveVariables}). This information changes the
probability
measure $\mu_0$ into a conditional measure $\mu_V$ with density
$f_V$
given by (\ref{InitialDensity}). Conditional averages of operators
$O[u(\cdot)]$ can be calculated by using the following three
lemmas:

\begin{lemma}
\label{Lemma1}
The conditional expectation of the function $u(x)$ is a linear form
in
the conditioning data $V$:
\begin{equation}
\average{u(x)}_V = \average{u(x)} +
\sum_{\alpha=1}^N c_\alpha(x) \BRK{V_\alpha -
\average{U_\alpha[u(\cdot)]}},
\label{ConditionalMean}
\end{equation}
where the vector of functions $c_\alpha(x)$ is given by
\begin{equation}
c_\alpha(x) = \sum_{\beta=1}^N
\scalar{a^{-1}(x,\cdot)}{g_\beta(\cdot)}
m^{-1}_{\beta\alpha},
\label{DefinitionC}
\end{equation}
and where the $m^{-1}_{\beta\alpha}$ are the entries of an $N\times
N$
matrix $M^{-1}$ whose inverse $M$ has entries
\begin{equation}
m_{\beta\alpha} =
\cov{U_\beta[u(\cdot)]}{U_\alpha[u(\cdot)]} =
\iint g_\beta(x) a^{-1}(x,y) g_\alpha(y) \,dx \, dy.
\end{equation}
\end{lemma}

{\bf Proof.}
Given the prior measure $\mu_0$ and the values $V$ of the
collective
variables $U$, we define a \emph{regression function} (an
approximant
to $u(x)$) of the form
\begin{equation}
R(x) = \sum_{\alpha=1}^N r_\alpha(x) V_\alpha + s(x),
\label{RegressionCurve}
\end{equation}
where the functions $r_\alpha(x)$ and $s(x)$ are chosen such to
minimize the mean square error,
\begin{equation}
E(x) = \average{e^2(x)} \equiv
\average{ \Brk{
u(x) - \sum_{\alpha=1}^N r_\alpha(x) U_\alpha[u(\cdot)] -
s(x)}^2}.
\end{equation}
for all $x$. Note that this is an \emph{unconditional} average with
respect
to $\mu_0$.

Minimization with respect to $s(x)$ implies that
\begin{equation}
\frac{\partial E(x)}{\partial s(x)} =
\average{e(x)} =
\average{ u(x) - \sum_{\alpha=1}^N r_\alpha(x) U_\alpha[u(\cdot)] -
s(x) } = 0,
\label{minimumB}
\end{equation}
which, combined with (\ref{RegressionCurve}), yields
\begin{equation}
R(x) = \average{u(x)} +
\sum_{\alpha=1}^N r_\alpha(x) \BRK{\average{U_\alpha[u(\cdot)]} -
V_\alpha}.
\label{ExpressionForR}
\end{equation}
Minimization with respect to $r_\alpha(x)$ implies:
\begin{equation}
\frac{\partial E(x)}{\partial r_\alpha(x)} =
\average{e(x) \, U_\alpha[u(\cdot)]} =
\average{
\Brk{u(x) - \sum_{\beta=1}^N r_\beta(x) U_\beta[u(\cdot)] - s(x)}
U_\alpha[u(\cdot)] } = 0.
\label{minimumA}
\end{equation}
Equation (\ref{minimumA}) can be rearranged by substituting
equations
(\ref{UncCovariance}) and (\ref{minimumB}) into it, and using the
fact
that $U_\alpha[u(\cdot)] = \scalar{g_\alpha(\cdot)}{u(\cdot)}$:
\begin{equation}
\sum_{\beta=1}^N \cov{U_\alpha[u(\cdot)]}{U_\beta[u(\cdot)]} r_\beta(x) =
\scalar{g_\alpha(\cdot)}{a^{-1}(x,.)}.
\end{equation}
One readily identifies the functions $r_\alpha(x)$ as satisfying
the
definition (\ref{DefinitionC}) of the functions $c_\alpha(x)$.
Comparing (\ref{ExpressionForR}) with (\ref{ConditionalMean}), the
regression function is nothing but the right-hand side of equation
(\ref{ConditionalMean}).

It remains to show that the regression curve equals also the
left-hand
side of (\ref{ConditionalMean}). Consider equation
(\ref{minimumA}):
it asserts that the random variable $e(x)$ is statistically
orthogonal
to the random variables $U_\alpha[u(\cdot)]$. Note that both
$e(x)$
and the collective variables $U_\alpha$ are linear functionals of
the
Gaussian function $u(x)$, and are therefore jointly Gaussian.
Jointly
Gaussian variables that are statistically orthogonal are
independent,
hence, the knowledge of the value assumed by the variables
$U_\alpha[u(\cdot)]$ does not affect the expectation value of
$e(x)$,
\begin{equation}
\average{u(x) - \sum_{\alpha=1}^N r_\alpha(x) U_\alpha[u(\cdot)] -
s(x)}_V =
\average{u(x) - \sum_{\alpha=1}^N r_\alpha(x) U_\alpha[u(\cdot)] -
s(x)}.
\end{equation}
The function $s(x)$ is not random and
$\average{U_\alpha[u(\cdot)]}_V
= V_\alpha$, from which immediately follows that
\begin{equation}
\average{u(x)}_V = \average{u(x)} +
\sum_{\alpha=1}^N r_\alpha(x) \BRK{V_\alpha -
\average{U_\alpha[u(\cdot)]}},
\end{equation}
This completes the proof.
\bigskip

\begin{lemma}
\label{Lemma2}
The conditional covariance of the function $u(x)$ differs from
the unconditional covariance by a function that depends on the
kernels
$g_\alpha(x)$, without reference to the conditioning data $V$:
\begin{equation}
\cov{u(x)}{u(y)}_V = \cov{u(x)}{u(y)} -
\sum_{\alpha=1}^N c_\alpha(x)
\scalar{g_\alpha(\cdot)}{a^{-1}(\cdot,y)}.
\label{ConditionalCovariance}
\end{equation}
\end{lemma}

{\bf Proof.}
The proof follows the same line as the second part of the proof of
Lemma~\ref{Lemma1}. Consider the following expression:
\begin{equation}
e(x) e(y) =
\Brk{u(x) - \sum_{\alpha=1}^N r_\alpha(x) U_\alpha[u(\cdot)] -
s(x)}
\Brk{u(y) - \sum_{\beta=1}^N  r_\beta(y)
U_\beta[u(\cdot)]  - s(y)}.
\end{equation}
Both $e(x)$ and $e(y)$ are independent of the collective variables
$U$. It is always true that if $A_1$, $A_2$ and $A_3$ are random
variables with $A_3$ being independent of $A_1$ and $A_2$, then
$\average{A_1 A_2}_{A_3} = \average{A_1 A_2}$. Hence,
\begin{equation}
\average{e(x) e(y)}_V = \average{e(x) e(y)},
\end{equation}
from which (\ref{ConditionalCovariance}) follows after
straightforward
algebra.
\bigskip

\begin{lemma}
\label{Lemma3}
Wick's theorem extends to conditional expectations:
\begin{multline}
\average{(u_{i_1}-\average{u_{i_1}}_V) \cdots
(u_{i_l}-\average{u_{i_l}}_V)}_V = \\
\left\{
\begin{aligned}
0, & \qquad l \text{ odd} \\
\sum
\cov{u_{i_{p_1}}}{u_{i_{p_2}}}_V \cdots
\cov{u_{i_{p_{l-1}}}}{u_{i_{p_l}}}_V, & \qquad l \text{ even}
\end{aligned}
\right. ,
\label{ConditionalWick}
\end{multline}
where again the summation is over all possible pairings of
$\{i_1,\ldots,i_l\}$.
\end{lemma}

{\bf Proof.}
Using the fact that a delta function can be represented as the
limit
of a narrow Gaussian function, the conditional expectation of any
list
of observables, $O_1[u(\cdot)],\ldots,O_p[u(\cdot)]$, can be
expressed
as
\begin{equation}
\average{O_1[u(\cdot)] \cdots O_p[u(\cdot)]}_V = \lim_{\Delta\to0}
\int O_1[u(\cdot)] \cdots O_p[u(\cdot)]\,f_V^\Delta[u(\cdot)]
[du],
\end{equation}
where
\begin{equation}
f_V^\Delta[u(\cdot)] =
c_\Delta \,f_0[u(\cdot)] \, \prod_{\alpha=1}^N
\frac{1}{\sqrt{\pi} \Delta} \exp
\Brk{ - \frac{\brk{U_\alpha[u(\cdot)] - V_\alpha}^2}{\Delta^2}},
\label{f_delta}
\end{equation}
the coefficient $c_\Delta$ is a normalization, and the order of
the
limit $\Delta\to0$ and the functional integration has been
interchanged. Note that the exponential in (\ref{f_delta}) is
quadratic in $u(x)$, hence the finite-$\Delta$ probability density
$f_V^\Delta[u(\cdot)]$ is Gaussian, Wick's theorem applies, and
the
limit $\Delta\to0$ can finally be taken.
\bigskip

The conditional expectation of any observable $O[u(\cdot)]$ can be
deduced, in principle, from a combination of Lemmas
\ref{Lemma1}-\ref{Lemma3}.

In the examples considered below, the dependent variable $u(x,t)$ is
a
vector; let $u^i(x,t)$ denote the $i$'th component of the
$d$-dimensional vector $u(x,t)$.  All the above relations are
easily
generalized to the vector case. To keep notations as clear as
possible, we denote indices associated with the collective
variables
by Greek subscripts, and indices associated with the components of
$u$
by Roman superscripts.  The probability density $f_0[u(\cdot)]$
is
Gaussian if it is of the following form,
\begin{equation}
f_0[u(\cdot)] = \frac{1}{Z} \exp\brk{
- \half \sum_{i,j=1}^d \iint u^i(x) a^{ij}(x,y) u^j(y)\,dx\,dy
+ \sum_{i=1}^d \int b^i(x) u^i(x) \, dx},
\label{FormOfF}
\end{equation}
where $a^{ij}(x,y)$ are now the entries of a $d \times d$ matrix
of
functions, and $b^i(x)$ are the entries of a vector of
functions. These functions are related to the mean and the
covariance
of the vector $u(x)$ by
\begin{equation}
\average{u^i(x)} = \sum_{j=1}^d
\scalar{[a^{-1}(x,\cdot)]^{ij}}{b^j(\cdot)},
\end{equation}
and
\begin{equation}
\cov{u^i(x)}{u^j(y)} = [a^{-1}(x,y)]^{ij},
\end{equation}
where $[a^{-1}(x,y)]^{ij}$ is defined by
\begin{equation}
\sum_{j=1}^d \scalar{[a^{-1}(x,\cdot)]^{ij}}{a^{jk}(\cdot,y)} =
\delta(x-y)\,\delta_{ik}.
\end{equation}

Suppose now that a set of measurements reveals the values
$V^i_\alpha$
of a matrix of collective variables of the form,
\begin{equation}
U^i_\alpha[u(\cdot)] = \scalar{g_\alpha(\cdot)}{u^i(\cdot)},
\end{equation}
where $\alpha=1,\ldots,N$ and $i=1,\ldots,d$. The conditional
expectation
and covariance of $u^i(x)$ are given by straightforward
generalizations of Lemmas~\ref{Lemma1} and \ref{Lemma2}:
\begin{equation}
\average{u^i(x)}_V = \average{u^i(x)} +
\sum_{\alpha=1}^N \sum_{j=1}^d c^{ij}_\alpha(x)
\BRK{V^j_\alpha - \average{U^j_\alpha[u(\cdot)]}},
\label{VecConditionalAverage}
\end{equation}
and
\begin{equation}
\cov{u^i(x)}{u^j(y)}_V = \cov{u^i(x)}{u^j(y)} -
\sum_{\alpha=1}^N \sum_{k=1}^d
c^{ik}_\alpha(x)  \scalar{g_\alpha(\cdot)}{[a^{-1}(\cdot,y)]^{kj}}.
\label{VecConditionalCovariance}
\end{equation}
where
\begin{equation}
c^{ij}_\alpha(x) = \sum_{\beta=1}^N \sum_{k=1}^d
\scalar{[a^{-1}(x,\cdot)]^{ik} }{g_\beta(\cdot)}
[m^{-1}]_{\beta\alpha}^{kj},
\label{}
\end{equation}
and where the $[m^{-1}]_{\beta\alpha}^{ij}$ are the entries of an
$N\times N\times d\times d$ tensor $M^{-1}$ whose inverse $M$ has
entries
\begin{equation}
m_{\beta\alpha}^{ij} =
\iint g_\beta(x) [a^{-1}(x,y)]^{ij} g_\alpha(y) \,dx \, dy.
\end{equation}

\section{A linear Schr\"odinger equation}
\label{SecLin}

{\bfseries The equations of motion.}  The first example is a
linear
Schr\"odinger equation that we write as a pair of real
equations:
\begin{equation}
\begin{split}
p_t &= - q_{xx} + m_0^2 q \\
q_t &= + p_{xx} - m_0^2 p \\
\end{split} ,
\label{Lin:Equation}
\end{equation}
where $p(x,t)$ and $q(x,t)$ are defined on the domain $(0,2\pi]$,
$m_0$
is a constant, and periodic boundary conditions are assumed.
Equations
(\ref{Lin:Equation}) are the Hamilton equations of motion for the
Hamiltonian \cite{FH65},
\begin{equation}
H[p(\cdot),q(\cdot)] = \half \int_0^{2\pi}
\Brk{(p_x)^2 + (q_x)^2 + m_0^2(p^2 + q^2)} \, dx,
\label{Lin:Hamiltonian}
\end{equation}
with $p(x)$ and $q(x)$ being the canonically conjugate variables.

{\bfseries The prior measure.}  Equation (\ref{Lin:Equation})
preserves
any density that is a function of the Hamiltonian. We will assume
that
the prior measure is given by the canonical density,
\begin{equation}
f_0[p(\cdot),q(\cdot)] = \exp \BRK{ - H[p(\cdot),q(\cdot)]},
\label{Lin:Density}
\end{equation}
where the temperature has been chosen equal to one.

The measure defined by equation (\ref{Lin:Density}) is absolutely
continuous with respect to a Wiener measure \cite{McK95}, and its
samples are, with probability one, almost nowhere differentiable.
The
corresponding solutions of the equations of motion are weak and
hard
to approximate numerically.

The Hamiltonian (\ref{Lin:Hamiltonian}) is quadratic in $p$ and
$q$,
hence the probability density (\ref{Lin:Density}) is Gaussian. By
symmetry we see that the unconstrained means $\average{p(x)}$ and
$\average{q(x)}$ are zero.  To extract the matrix of
covariance
functions $A^{-1}$, we write the Hamiltonian
(\ref{Lin:Hamiltonian})
as a double integral:
\begin{multline}
H[p(\cdot),q(\cdot)] =
\iint
\bigg[ p_x(x) \delta(x-y) p_x(y) + q_x(x) \delta(x-y) q_x(y) + \\
+ m_0^2\, p(x) \delta(x-y) p(y) + m_0^2\, q(x) \delta(x-y) q(y)
\bigg]
\, dx\,dy.
\end{multline}
Integration by parts shows that the entries of the matrix of
functions
$A$ are
\begin{equation}
a^{ij}(x,y) = \Brk{ - \delta''(x-y) + m_0^2\,\delta(x-y)}
\,\delta_{ij},
\end{equation}
where the indices $i$ and $j$ represent either $p$ or $q$, and
$\delta''(\cdot)$ is a second derivative of a delta function. The
integral equation for the inverse operator $A^{-1}$ can be solved
by
Fourier series. The result is a translation-invariant diagonal
matrix
\begin{equation}
[a^{-1}(x,y)]^{ij} = \frac{1}{2\pi} \delta_{ij}
\sum_{k=-\infty}^\infty
\frac{e^{i k(x-y)}}{k^2 + m_0^2}.
\label{Lin:Covariance}
\end{equation}

{\bfseries The collective variables.}  We assume that the
initial data
for equations (\ref{Lin:Equation}) are drawn from the distribution
(\ref{Lin:Density}), and that $2N$ measurements have revealed the
values of the $2N$ collective variables,
\begin{equation}
\begin{split}
U^p_\alpha[p(\cdot),q(\cdot)] &\equiv
\scalar{g_\alpha(\cdot)}{p(\cdot)}
= V^p_\alpha \\
U^q_\alpha[p(\cdot),q(\cdot)] &\equiv
\scalar{g_\alpha(\cdot)}{q(\cdot)}
= V^q_\alpha \\
\end{split}   ,
\label{Lin:Collective}
\end{equation}
for $\alpha=1,\ldots,N$. The kernels $g_\alpha(x)$ are translates
of
each other, $g_\alpha(x) = g(x-x_\alpha)$, and the points
$x_\alpha=2\pi\alpha/N$ form a regular mesh on the interval
$(0,2\pi]$. We choose
\begin{equation}
g(x) = \frac{1}{\sqrt{\pi}\sigma} \sum_{\tau=-\infty}^\infty
\exp\Brk{-\frac{(x-2\pi\tau)^2}{\sigma^2}},
\end{equation}
i.e., the kernel is a normalized Gaussian whose width is $\sigma$,
with suitable images to enforce periodicity. The Fourier
representation of $g(x)$ is
\begin{equation}
g(x) = \frac{1}{2\pi} \sum_{k=-\infty}^\infty
e^{i k x} e^{-\frac{1}{4} k^2 \sigma^2}.
\label{FourierSeriesOfG}
\end{equation}

We could have trivialized this example by choosing as kernels a set
of
trigonometric functions, which are eigenfunctions of the evolution
operator. The goal here is to demonstrate what one could do when
an
exact representation of the eigenfunctions is not known.

{\bfseries Conditional expectation.}  We now demonstrate the
application of the Lemmas derived in the previous section. Given
the
initial data, $V^p$ and $V^q$, we may calculate the expectation of
the
functions $p(x)$ and $q(x)$; these conditional averages are given
by
equation (\ref{VecConditionalAverage}). Because the unconditional
averages of $p(x)$, $q(x)$, $U^p_\alpha$ and $U^q_\alpha$ all
vanish,
and the unconditional covariance $[a^{-1}(x,y)]^{ij}$ is diagonal
with
respect to $i$ and $j$ ($p$ and $q$ are independent), equation
(\ref{VecConditionalAverage}) reduces to a simpler expression; the
conditional average of $p(x)$, for example, is
\begin{equation}
\average{p(x)}_V = \sum_{\alpha=1}^N c^{pp}_\alpha(x) V^p_\alpha,
\label{Lin:Interpolation}
\end{equation}
where
\begin{equation}
c^{pp}_\alpha(x) = \sum_{\beta=1}^N
\scalar{[a^{-1}(x,\cdot)]^{pp}}{g_\beta(\cdot)} [m^{-1}]^{pp}_{\beta\alpha}
= c^{qq}_\alpha(x),
\end{equation}
and $[m^{-1}]^{pp}_{\beta\alpha}$ are the entries of an $N\times
N$
matrix $M^{-1}$ (the upper indices $p$ are considered as fixed)
whose
inverse $M$ has entries
\begin{equation}
m^{pp}_{\beta\alpha} = \iint
g_\beta(x) [a^{-1}(x,y)]^{pp} g_\alpha(y)\,dx\,dy =
m^{qq}_{\beta\alpha}.
\end{equation}
Substituting the Fourier representations of $A^{-1}$
(\ref{Lin:Covariance}) and $g$ (\ref{FourierSeriesOfG}), we obtain
\begin{equation}
c^{pp}_\alpha(x) = \frac{1}{2\pi} \sum_{\alpha=1}^N
\sum_{k=-\infty}^\infty \frac{e^{-\frac{1}{4} k^2 \sigma^2}}{k^2 +
m_0^2}
\exp\Brk{i k(x - x_\beta)} [m^{-1}]^{pp}_{\beta\alpha},
\end{equation}
and
\begin{equation}
m^{pp}_{\beta\alpha} = \frac{1}{2\pi}
\sum_{k=-\infty}^\infty \frac{e^{-\half k^2 \sigma^2}}{k^2 +
m_0^2}
\exp\Brk{i k(x_\alpha - x_\beta)}.
\end{equation}

The regression function (\ref{Lin:Interpolation}) can be viewed as
an
``optimal interpolant''; it is the expectation value of the
function
$p(x)$ given what is known. Examples of regression functions are
plotted in Figure~\ref{fig:LinInterp} for a mesh of $N=5$ points.
The
open circles represent the values of the five collective variables
$V^p_\alpha$; the abscissa is the location of the point $x_\alpha$
around which the average is computed, and the ordinate is the value
of
the corresponding collective variable. The three curves represent
the
interpolating function (\ref{Lin:Interpolation}) for three
different
values of the kernel width: $\sigma=\Delta x = 2\pi/N$ (solid
line),
$\sigma=0.5\,\Delta x$ (dashed line), and $\sigma=0.1\,\Delta x$
(dash-dot line).  The parameter $m_0$ was taken to be one.

\begin{figure}
\caption{Example of regression functions for the linear
Schr\"odinger
equation. Values for five collective variables were chosen,
representing local averages of $p(x)$ on a uniformly spaced grid.
The
kernels are translates of each other and have Gaussian profiles of
width $\sigma$ centered at the grid points. The lines represent
the
regression function, or optimal interpolant $\average{p(x)}_V$ given
by
equation (\ref{Lin:Interpolation}) for $\sigma=\Delta x$ (solid),
$\sigma=0.5\,\Delta x$ (dashed), and $\sigma=0.1\,\Delta x$
(dash-dot).}
\label{fig:LinInterp}
\end{figure}

{\bfseries Time evolution.} We next consider the time evolution of
the
mean value of the collective variables $U^p$ and $U^q$, first based
on
the approximating scheme (\ref{TheEffectiveEquation}).  The
equation
for $V_\alpha^p$, for example, is
\begin{equation}
\begin{split}
\frac{d V^p_\alpha}{d t} &=
\average{\scalar{g_\alpha(\cdot)}{-q_{xx}
(\cdot) + m_0^2 q(\cdot)}}_V =
\\
& = - \scalar{g_\alpha(\cdot)}{\frac{\partial^2}{\partial x^2}
\average{q(\cdot)}_V} + m_0^2
\scalar{g_\alpha(\cdot)}{\average{q(\cdot)}_V}.
\end{split}
\end{equation}
Substituting the regression function (\ref{Lin:Interpolation}) we
find:
\begin{equation}
\frac{d V^p_\alpha}{d t} = \sum_{\gamma=1}^N
\BRK{\sum_{\beta=1}^N
\scalar{g_\alpha(\cdot)}{g_\beta(\cdot)}
[m^{-1}]^{qq}_{\beta\gamma}}
V^q_{\gamma}.
\label{Lin:Effective}
\end{equation}
A similar equation is obtained for $V^q_\alpha$ by the symmetry
transformation $V^p_\alpha \to V^q_\alpha$ and $V^q_\alpha \to -
V^p_\alpha$. Equation (\ref{Lin:Effective}) represents a set of
$2N$
ordinary differential equations that approximate the mean evolution
of
the collective variables. These equations are easy to solve with
standard ODE solvers. Note that the matrix elements in braces need
to
be computed only once to define the scheme.

We next calculate the \emph{exact} mean value of the collective
variables, $U^p$ and $U^q$, at time $t$, conditioned by the
initial
data, $V^p$ and $V^q$, at time $t=0$, so that they can be compared
with the result $V(t)$ of the scheme we just presented. We are able
to
do so in the present case because the equations are linear, and a
simple representation of the evolution operator can be found.

The solution to the initial value problem (\ref{Lin:Equation}) can
be
represented by Fourier series,
\begin{equation}
\begin{split}
p(x,t) &= \frac{1}{2\pi} \sum_{k=-\infty}^\infty
\int e^{ik(x-y)} \Brk{ p(y)\,\cos \omega t +
q(y)\,\sin \omega t} \, dy \\
q(x,t) &= \frac{1}{2\pi} \sum_{k=-\infty}^\infty
\int e^{ik(x-y)} \Brk{ q(y)\,\cos \omega t -
p(y)\,\sin \omega t} \, dy \\
\end{split}
\label{Lin:Exact1}
\end{equation}
where $p(y)$ and $q(y)$ are the (random) initial conditions, and
$\omega=k^2 + m_0^2$.

The expectation values of the collective variables $U^p_\alpha$
and
$U^q_\alpha$ are obtained by averaging the scalar products
$\scalar{p(\cdot,t)}{g_\alpha(\cdot)}$ and
$\scalar{q(\cdot,t)}{g_\alpha(\cdot)}$ with respect to the initial
distribution. Because equations (\ref{Lin:Exact1}) are linear in
the
random variables $p(y)$ and $q(y)$ this gives
\begin{equation}
\begin{split}
\average{U^p_\alpha[p(\cdot),q(\cdot),t]}_V &=
\frac{1}{2\pi} \sum_{k=-\infty}^\infty
\int e^{ik(x_\alpha-y)-\frac{1}{4}k^2 \sigma^2}
\Brk{\average{p(y)}_V \cos \omega t + \average{q(y)}_V \sin \omega t}
dy \\
\average{U^q_\alpha[p(\cdot),q(\cdot),t]}_V & =
\frac{1}{2\pi} \sum_{k=-\infty}^\infty
\int e^{ik(x_\alpha-y)-\frac{1}{4}k^2 \sigma^2}
\Brk{\average{q(y)}_V \cos \omega t - \average{p(y)}_V \sin \omega t}
dy \\
\end{split} .
\label{Lin:Exact2}
\end{equation}
Note that in the linear case averaging and time evolution commute;
equation (\ref{Lin:Exact2}) would have also been obtained if we
first
computed the mean initial state, $\average{p(y)}_V$ and
$\average{q(y)}_V$, evolved it in time according to
(\ref{Lin:Exact1}), and finally computed the collective variables
by
taking the appropriate scalar products.

To complete the calculation, we substitute the linear regression
formula (\ref{Lin:Interpolation}) for $\average{p(y)}_V$ and
$\average{q(y)}_V$ and obtain:
\begin{equation}
\begin{split}
\average{U^p_\alpha[p(\cdot),q(\cdot),t]}_V &=
\sum_{\beta,\gamma=1}^N \BRK{
c^C_{\alpha\beta}(t) [m^{-1}]^{pp}_{\beta\gamma} V^p_\gamma +
c^S_{\alpha\beta}(t) [m^{-1}]^{qq}_{\beta\gamma} V^q_\gamma} \\
\average{U^q_\alpha[p(\cdot),q(\cdot),t]}_V &=
\sum_{\beta,\gamma=1}^N \BRK{
c^C_{\alpha\beta}(t) [m^{-1}]^{pp}_{\beta\gamma} V^q_\gamma -
c^S_{\alpha\beta}(t) [m^{-1}]^{qq}_{\beta\gamma} V^p_\gamma} \\
\end{split} ,
\label{Lin:Exact3}
\end{equation}
where
\begin{equation}
c^C_{\alpha\beta}(t) = \frac{1}{2\pi} \sum_{k=-\infty}^\infty
\frac{\cos \omega t}{\omega} e^{ik(x_\alpha-x_\beta)}
e^{-\half k^2 \sigma^2},
\end{equation}
and
\begin{equation}
c^S_{\alpha\beta}(t) = \frac{1}{2\pi} \sum_{k=-\infty}^\infty
\frac{\sin \omega t}{\omega} e^{ik(x_\alpha-x_\beta)}
e^{-\half k^2 \sigma^2}.
\end{equation}

{\bfseries Results.} We now compare the exact formula
(\ref{Lin:Exact3}) for the future expectation value of the
collective
variables to the approximation (\ref{Lin:Effective}). Figures
\ref{Fig:LinEvolution}$a$--\ref{Fig:LinEvolution}$c$ compare between the two
evolutions for $N=5$ and randomly selected initial data,
$V^p_\alpha$
and $V^q_\alpha$. The graphs show the mean time evolution of the
collective variable $U^p_1[p(\cdot),q(\cdot)]$.  The same set
of
initial values was used in the three plots; the difference is in
the
width $\sigma$ of the kernels $g_\alpha(x)$: $\sigma=\Delta x$
(Figure
\ref{Fig:LinEvolution}$a$), $\sigma=0.5\,\Delta x$ (Figure
\ref{Fig:LinEvolution}$b$), and $\sigma=0.1\,\Delta x$ (Figure
\ref{Fig:LinEvolution}$c$). In the first case, in which the kernel
width
equals the grid spacing, the approximation is not distinguishable
from
the exact solution on the scale of the plot for the duration of
the
calculation. The two other cases show that the narrower the kernel
is,
the sooner the curve deviates from the exact solution.

\begin{figure}
\caption{Mean evolution of the collective variable
$U^p_1[p(\cdot),q(\cdot)]$ for $N=5$, and a random choice of the
initial data $V^p$ and $V^q$. The open dots represent the exact
solution (\ref{Lin:Exact3}), whereas the lines represent the
approximate solution obtained by an integration of the set of $10$
ordinary differential equations (\ref{Lin:Effective}). The three
graphs are for different values of the kernel width $\sigma$: (a)
$\sigma=\Delta x$, (b) $\sigma=0.5\,\Delta x$, and (c)
$\sigma=0.1\,\Delta x$.}
\label{Fig:LinEvolution}
\end{figure}

\section{A nonlinear Hamiltonian system}
\label{SecNonlin}

{\bfseries The equations of motion.}  The method demonstrated in
the
preceding section can be generalized to a nonlinear
Schr\"odinger
equation. However, we want to exhibit the power of our method by
comparing the solutions that it yields to exact solutions; in the
nonlinear case, exact solutions of problems with random initial
conditions are hard to find, so we resort to a stratagem. Even
though
our method applies to nonlinear partial differential equations, we
study instead a finite dimensional system of $2n$ ordinary
differential equations that is formally a finite difference
approximation of a nonlinear Schr\"odinger equation:
\begin{equation}
\begin{split}
\frac{d p(j)}{d t} =  -\frac{q(j-1) - 2q(j) + q(j+1)}{\Delta
x^2}
+ q^3(j) \\
\frac{d q(j)}{d t} =  +\frac{p(j-1) - 2p(j) + p(j+1)}{\Delta
x^2}
- p^3(j) \\
\end{split}
\qquad j=1,\ldots,n,
\label{NSE:Equation}
\end{equation}
where $\Delta x=1/n$ is the mesh spacing, and periodicity is
enforced
with $p(0) \equiv p(n)$, $p(n+1) \equiv p(1)$, etc; this system is
non-integrable for $n,1$.  The approximation is only formal
because we
shall be considering non-smooth data which give rise to weak
solutions
that cannot be readily computed by difference methods.

We shall pretend that $n$ is so large that the system
(\ref{NSE:Equation}) cannot be solved on a computer, and shall
therefore seek an approximation that requires a computation with
fewer
variables. In practice we shall pick an $n$ small enough so that
the
results of the approximate procedure can be compared to an
\emph{ensemble} of exact solutions.

{\bfseries The prior measure.} The system of equations
(\ref{NSE:Equation}) is the Hamilton equations of motion for the
Hamiltonian
\begin{equation}
H[p,q] = \half \sum_{j=1}^n \BRK{
\Brk{\frac{p(j+1)-p(j)}{\Delta x}}^2 +
\Brk{\frac{q(j+1)-q(j)}{\Delta x}}^2 +
\half \Brk{p^4(j) + q^4(j)}},
\label{NSE:Hamiltonian}
\end{equation}
where $p \equiv(p(1),\ldots,p(n))$ and $q
\equiv(q(1),\ldots,q(n))$. The differential equations
(\ref{NSE:Equation}) preserve the canonical density
\begin{equation}
f_0[p,q] = \exp \BRK{-H[p,q]},
\label{NSE:Prior}
\end{equation}
which we postulate, as before, to be the prior probability
density.

The prior density (\ref{NSE:Prior}) is not Gaussian, which raises
a
technical difficulty in computing expectation values. We adopt here
an
approximate procedure where the density (\ref{NSE:Prior}) is
approximated by a Gaussian density that yields the same first and
second moments (means and covariances) of the vectors $p$ and $q$.
The
means are zero by symmetry:
\begin{equation}
\average{p(j)} = \average{q(j)} = 0
\end{equation}
(positive and negative values of these have equal weight). Also
all
$p$'s and $q$'s are uncorrelated:
\begin{equation}
\average{p(j_1) q(j_2)} = 0,
\end{equation}
since the density factors into a product of a density for the
$p$'s
and a density for the $q$'s. Thus $\average{p(j_1) p(j_2)}
=\average{q(j_1) q(j_2)}$ are  the only non-trivial
covariances. Finally, since the Hamiltonian is translation
invariant,
these covariances depend only on the separation between the indices
$j_1$
and $j_2$, and are symmetric in $j_1-j_2$.

To relate the present discrete problem to the continuous formalism
used in the preceding section we write in analogy to
(\ref{Lin:Covariance})
\begin{equation}
\begin{split}
\cov{p(j_1)}{p(j_2)} &= [a^{-1}(j_1,j_2)]^{pp} = c(|j_1-j_2|)
\\
\cov{p(j_1)}{q(j_2)} &= [a^{-1}(j_1,j_2)]^{pq} = 0, \\
\label{NSE:Covariance}
\end{split}
\end{equation}
with $j_1,j_2 = 1,\ldots,n$.  We computed the numbers,
$c(|j_1-j_2|)$,
for $n=16$ and $j_1-j_2=0,\ldots,15$ by a Metropolis Monte-Carlo
algorithm \cite{BH92}; the covariances obtained this way are shown
in
Figure \ref{Fig:NSECovariance}. Along with the zero means, the
numbers
represented in Figure \ref{Fig:NSECovariance} completely determine
the
\emph{approximate} prior distribution.

\begin{figure}
\caption{The covariance $\average{p(i)\,p(j)} = \average{q(i)\,q(j)}$
as
function of the grid separation $i-j$ for the non-Gaussian
probability
distribution (\ref{NSE:Prior}) with $n=16$. These values were
computed
by a Metropolis Monte-Carlo simulation.}
\label{Fig:NSECovariance}
\end{figure}

{\bfseries The collective variables.} We next define a set of $2N$
collective variables ($N$), whose values we assume to be given
at
the initial time. The class of collective variables that is the
discrete analog of (\ref{Lin:Collective}) is of the form
\begin{equation}
\begin{split}
U^p_\alpha[p,q] = \scalar{g_\alpha(\cdot)}{p(\cdot)} \equiv
\sum_{j=1}^n g_\alpha(j) p(j) \\
U^q_\alpha[p,q] = \scalar{g_\alpha(\cdot)}{q(\cdot)} \equiv
\sum_{j=1}^n g_\alpha(j) q(j) \\
\end{split}
\qquad \alpha=1,\ldots,N,
\label{NSE:Collective}
\end{equation}
where the $g$'s are discrete kernels. In the calculations we
exhibit
we chose $n=16$ and $N=2$ so that we aim to reduce the number of
degrees of freedom by a factor of $8$. We pick as kernels
discretized
Gaussian functions centered at the grid points $j=1$ and $j=9$:
\begin{equation}
\begin{split}
g_1(j) = \frac{1}{Z} \exp\BRK{-
\frac{d^2(1,j)}{n^2 \sigma^2}} \\
g_2(j) = \frac{1}{Z} \exp\BRK{-
\frac{d^2(9,j)}{n^2 \sigma^2}} \\
\end{split}
\label{NSE:Kernels}
\end{equation}
where $Z$ is a normalizing constant, $\sigma=0.25$, and
$d(j_1,j_2)$
is a distance function over the periodic index axis, i.e., it is
the
minimum of $|j_1-j_2|$, $|j_1-j_2-n|$, and $|j_1-j_2+n|$.

{\bfseries Conditional expectation.} With the approximate measure
defined by the covariances (\ref{NSE:Covariance}), and the
collective
variables (\ref{NSE:Collective}), whose measured values are again
denoted by $V^p_\alpha$ and $V^q_\alpha$, we can approximate the
conditional expectation of various observables $O[p,q]$. We shall
need
specifically the conditional expectation values of $p(j)$ and
$p^3(j)$.

The approximate conditional expectation value of $p(j)$ is given
by
the discrete analog of equation (\ref{Lin:Interpolation}), namely,
\begin{equation}
\average{p(j)}_V = \sum_{\alpha=1}^N c^{pp}_\alpha(j) V^p_\alpha,
\label{NSE:cond1}
\end{equation}
where
\begin{equation}
c^{pp}_\alpha(j) = \sum_{\beta=1}^N
\scalar{[a^{-1}(j,\cdot)]^{pp}}{g_\beta(\cdot)}
[m^{-1}]^{pp}_{\beta\alpha},
\end{equation}
and
\begin{equation}
m^{pp}_{\beta\alpha} = \sum_{j_1,j_2=1}^n
g_\beta(j_1) [a^{-1}(j_1,j_2)]^{pp} g_\alpha(j_2).
\end{equation}
(Again, the matrix inversion is only with respect to the lower
indices
$\alpha$ and $\beta$.)

To calculate the approximate conditional expectation value of
$p^3(j)$
we first use Wick's theorem (Lemma \ref{Lemma3}):
\begin{equation}
\average{p^3(j)}_V = 3 \average{p^2(j)}_V \average{p(j)}_V -
2 \average{p(j)}_V^3,
\label{NSE:cond2}
\end{equation}
and then calculate the conditional second moment by using the
discrete
analog of equation (\ref{ConditionalCovariance}):
\begin{equation}
\average{p^2(j)}_V = \average{p(j)}_V^2 +
[a^{-1}(j,j)]^{pp} - \sum_{\alpha=1}^N
c^{pp}_\alpha(j)
\scalar{g_\alpha(\cdot)}{[a^{-1}(\cdot,j)]^{pp}}.
\end{equation}

{\bfseries Time evolution.} The approximating scheme for
calculating
the mean evolution of the $2N$ collective variables $U^p$ and $U^q$
is
derived by substituting the kernels (\ref{NSE:Kernels}) and the
equations of motion (\ref{NSE:Equation}) in the approximation
formula
(\ref{TheEffectiveEquation}). The equation for $V^p_\alpha$, for
example, is
\begin{equation}
\begin{split}
\frac{d V^p_\alpha}{d t} &=
-\frac{1}{\Delta x^2} \sum_{j=1}^n g_\alpha(j)
\Brk{\average{q(j-1)}_V - 2 \average{q(j)}_V + \average{q(j+1)}_V}
+ \\ &+
\sum_{j=1}^n g_\alpha(j) \average{q^3(j)}_V.
\end{split}
\end{equation}
Substituting the expressions for the conditional expectations
(\ref{NSE:cond1}) and (\ref{NSE:cond2}), and performing the
summation,
using the values of the covariances plotted in
Figure~\ref{Fig:NSECovariance}, we explicitly obtain a closed set
of
$4$ ordinary differential equations. The equation for $V^p_1$ is:
\begin{equation}
\begin{split}
\frac{d V^p_1}{d t} &=
-19.5 \brk{V^q_2 - V^q_1} + \\ &+
\Brk{
1.50\, (V^q_1)^3 - 0.88\, (V^q_1)^2 V^q_2 + 0.27\,V^q_1 (V^q_2)^2
+
0.11\,(V^q_2)^3 }.
\end{split}
\label{NSE:Approx}
\end{equation}
The equation for $V^p_2$ is obtained by substituting
$1\leftrightarrow2$; the equations for $V^q_1$ and $V^q_2$ are
obtained by the transformation $p\to q$ and $q\to -p$.

Unlike in the linear case, we cannot calculate analytically the
mean
evolution of the collective variables. To assess the accuracy of
the
approximate equation (\ref{NSE:Approx}) we must compare the
solution
it yields with an average over an ensemble of solutions of the
``fine
scale'' problem (\ref{NSE:Equation}). To this end, we generated a
large number of initial conditions that are consistent with the
given
values, $V^p$ and $V^q$, of the collective variables. The
construction
of this ensemble was done by a Metropolis Monte Carlo algorithm,
where
new states are generated randomly by incremental changes, and
accepted
or rejected with a probability that ensures that for large enough
samples the distribution converges to the conditioned canonical
distribution. We generated an ensemble of $10^4$ initial
conditions;
each initial state was then evolved in time using a fourth-order
Runge-Kutta method. Finally, for each time level we computed the
distribution of collective variables, $U^p$ and $U^q$; the average
of
this distribution should be compared with the prediction of
equations
(\ref{NSE:Approx}).

The comparison between the true and the approximate evolution is
shown
in Figure~\ref{Fig:NSEEvolution}. Once again the reduced system of
equations reproduces the average behavior of the collective
variables
with excellent accuracy, but at a very much smaller computational
cost. Indeed, we compare one solution of $4$ equations to $10^4$
solutions of $32$ equations.

\begin{figure}
\caption{Evolution in time of the mean value of the four
collective
variable: $V^p_1$ ($\blacktriangledown$), $V^p_2$
($\blacktriangle$),
$V^q_1$ ($\blacksquare$), and $V^q_2$ ($\blacklozenge$). The
symbols
represent the values of these quantities obtained by solving the
$32$
equations (\ref{NSE:Equation}) for $10^4$ initial conditions
compatible with the initial data, and averaging. The solid lines
are
the values of the four corresponding functions obtained by
integrating
equation (\ref{NSE:Approx}). Figures (a) and (b) are for the time
intervals $[0,1]$ and $[0,10]$ respectively.}
\label{Fig:NSEEvolution}
\end{figure}

In Figure~\ref{Fig:NSESpread} we show the evolution of the
\emph{distribution} of values assumed by the collective variable
$U^p_1$; the data was extracted from the evolution of the
ensemble.
The distribution is initially sharply peaked, and spreads out as
time
evolves; yet, it remains sufficiently narrow throughout this
computation, so that the approximation that projects that
distribution
back onto a sharp one is reasonable. This indicates that the choice
of
collective variables, or kernels, was appropriate. The use of
narrow
kernels, or even point values, would have yielded a distribution
of
value that spreads out almost instantaneously.

\begin{figure}
\caption{Evolution of the distribution of the collective variable
$U^p_1$. The $x$-axis represents time, the $y$-axis represents the
value of $U^p_1$, and the $z$-axis is proportional to the density
of
states that correspond to the same value of $U^p_1$ at the given
time.}
\label{Fig:NSESpread}
\end{figure}

\section{Conclusions}

We have shown how to calculate efficiently, for a class of
problems,
the average behavior of an ensemble of solutions the individual
members of which are very difficult to evaluate. The approach is
reminiscent of statistical mechanics, where it is often easier to
predict the evolution of a mole of particles than to predict the
evolution of, say, a hundred particles, if one is content with the
average behavior of a set of coarse variables (collective
variables). The key step is the identification of a correspondence
between underresolution and statistics; underresolved data define,
together with prior statistical information, an ensemble of
initial
conditions, and the most one can aim for is to predict the
expectation
with respect to this ensemble of certain observables at future
times.
Our approach applies in those cases where prior statistical
information is available, and is consistent with the differential
equations; for example, it may consist of a measure invariant
under
the flow defined by the differential equations.  Fortunately,
there
are important classes of problems where we can find such
information.

We proposed a scheme (\ref{TheEffectiveEquation}) that advances in
time a set of variables that approximate the expectation values of
a
set of collective variables. As we explained, this scheme has to
be
viewed as a first approximation; more sophisticated schemes may be
designed by allowing the kernels to vary in time and/or by keeping
track of higher moments of the collective variables. Such
refinements
are the subject of ongoing research \cite{CKKT}.

One limitation of our present scheme can be perceived by
considering
the long time behavior of the nonlinear Hamiltonian system
presented
in Section~\ref{SecNonlin}. The flow induced by equations
(\ref{NSE:Equation}) is likely to be ergodic, hence the
probability
density function will approach, as $t\to\infty$, the invariant
distribution. Indeed, the initial data have a decreasing influence
on
the statistics of the solutions as time progresses. This implies
that
the expectation values of the observables $U^p$ and $U^q$ will tend
to
their unconditional means, i.e., will decay to zero. On the other
hand, no such decay occurs if one integrates the effective
equations
(\ref{NSE:Approx}) for very long times. One must conclude that the
present model is accurate for time intervals that are not longer
than
the time during which the initial data influence the outcome of
the
calculation.

The above discussion raises a number of questions interesting on
their
own: What is the range of influence, or the predictive power, of a
given set of data? How much information is contained in partial
data?
These questions need to be formulated in a more quantitative way;
they
are intimately related to the question of how to choose
appropriate
collective variables, and their scope is beyond any particular
method
of solution.

Finally, a full knowledge of the prior measure is a luxury one
cannot
always expect. One needs to consider problems where the
statistical
information is only partial; for example, a number of moments may
be
known from asymptotics and scaling analyses (e.g., in turbulence
theory \cite{Bar96,BC97,BC98a}). One can readily see from the
nonlinear example that one can make do with the knowledge of
means,
covariances, and perhaps some higher-order moments. In addition,
this
knowledge is needed only on scales comparable with the widths of
the
kernels.


\providecommand{\bysame}{\leavevmode\hbox to3em{\hrulefill}\thinspace}

\end{document}